%% file: agt-1-20.tex
\documentclass{gtart}
\input agtout

\lognumber{20}
\volumenumber{1}
\volumeyear{2001}
\papernumber{20}
\pagenumbers{411}{426}
\received{27 March 2001}
\accepted{6 July 2001}
\published{9 July 2001}

\usepackage{amsmath,amssymb,amscd}

\newtheorem{theorem}{Theorem}[section]
\newtheorem{lemma}[theorem]{Lemma}
\newtheorem{proposition}[theorem]{Proposition}
\newtheorem{corollary}[theorem]{Corollary}

\theoremstyle{definition}

\newtheorem{example}{Example}[section]

\newcommand\bfa{{\bf a}}
\newcommand\bfx{{\bf x}}
\newcommand\Z{{\mathbb Z}}
\newcommand\Q{{\mathbb Q}}

\newcommand\chiorb{\chi^{\text{orb}}}
\newcommand\degree{\operatorname{deg}}
\newcommand\genus{\operatorname{genus}}

\begin{document}
\title[$\pi_1$-injective surfaces in graph manifolds]
{Immersed and virtually embedded $\pi_1$-injective\\surfaces in
  graph manifolds}
\asciititle{Immersed and virtually embedded pi_1-injective surfaces in
  graph manifolds}

\author{Walter D. Neumann}
\address{Department of Mathematics, Barnard College, Columbia 
University\\New York, NY 10027, USA}
\email{neumann@math.columbia.edu}
\asciikeywords{pi_1-injective surface, graph manifold, separable,
  surface subgroup}

\begin{abstract}
  We show that many 3-manifold groups have no nonabelian surface
  subgroups. For example, any link of an isolated complex surface
  singularity has this property. In fact, we determine the exact class
  of closed graph-manifolds which have no immersed $\pi_1$-injective
  surface of negative Euler characteristic.  We also determine the
  class of closed graph manifolds which have no finite cover
  containing an embedded such surface.  This is a larger class.  Thus,
  manifolds $M^3$ exist which have immersed $\pi_1$-injective surfaces
  of negative Euler characteristic, but no such surface is virtually
  embedded (finitely covered by an embedded surface in some finite
  cover of $M^3$).
\end{abstract}
\asciiabstract{
  We show that many 3-manifold groups have no nonabelian surface
  subgroups. For example, any link of an isolated complex surface
  singularity has this property. In fact, we determine the exact class
  of closed graph-manifolds which have no immersed pi_1-injective
  surface of negative Euler characteristic.  We also determine the
  class of closed graph manifolds which have no finite cover
  containing an embedded such surface.  This is a larger class.  Thus,
  manifolds M^3 exist which have immersed pi_1-injective surfaces
  of negative Euler characteristic, but no such surface is virtually
  embedded (finitely covered by an embedded surface in some finite
  cover of M^3).}
\maketitle
\section{Introduction}
\label{sec:introduction}

It is widely expected that any closed $3$-manifold $M^3$ with infinite
fundamental group contains immersed $\pi_1$-injective surfaces. In
fact, standard conjectures of Waldhausen, Thurston and others imply
that some finite cover of $M^3$ has embedded $\pi_1$-injective
surfaces. If $M^3$ is hyperbolic --- or just simple and
non-Seifert-fibered, i.e., conjecturally hyperbolic by the
Geometrization Conjecture --- then an immersed $\pi_1$-injective
surface must have negative Euler characteristic.

We show here that many 3-manifolds have no immersed $\pi_1$-injective
surfaces of negative Euler characteristic and that yet more
3-manifolds have no virtually embedded ones (an immersion of a surface
$S$ in $M^3$ is a \emph{virtual embedding} if it can be lifted to an
embedding of a finite cover of $S$ in some finite cover of $M^3$).
Minimal surface theory implies that any $\pi_1$-injective surface in
an irreducible $3$-manifold is homotopic to an immersed one.

Our results might suggest caution for the standard conjectures about
hyperbolic manifolds. But the manifolds we study are graph-manifolds,
that is, $3$-manifolds with no non-Seifert-fibered pieces in their
JSJ-decomposition.  Thus, our results probably emphasise the very
different behaviour of hyperbolic manifolds and graph manifolds rather
than suggesting anything about what happens for hyperbolic $M^3$.

It is known that an immersed $\pi_1$-injective surface $S$ of
non-negative Euler characteristic is, up to isotopy, a collection of
tori and Klein bottles immersed parallel to fibers in Seifert fibered
pieces of the JSJ decomposition of $M^3$.  It follows with little
difficulty that $S$ is virtually embedded.  Thus an immersed
$\pi_1$-injective surface which is not virtually embedded must have
negative Euler characteristic.

The fact that graph manifolds can contain immersed $\pi_1$-injective
surfaces which are not virtually separable (homotopic to virtually
embedded) was first shown by H. Rubinstein and S. Wang
\cite{rubinstein-wang}, who in fact give a simple necessary and
sufficient criterion for a given horizontal immersed surface to be
virtually separable (the surface is \emph{horizontal} if it is
transverse to the fibers of the Seifert fibered pieces of $M^3$; this
implies $\pi_1$-injective).  They also show that, if a horizontal
surface is virtually separable in $M^3$, then it is separable: it
itself (rather than just some finite covering of it) lifts up to
homotopy to an embedding in some finite cover of $M^3$.  An embedded
horizontal surface in a graph manifold $M^3$ is a fiber of a fibration
of $M^3$ over $S^1$.  A necessary and sufficient condition for virtual
fibration of a graph manifold was given in \cite{neumann}.

Any infinite surface subgroup of $\pi_1(M^3)$ comes from a
$\pi_1$-injective immersion of a surface to $M^3$.  Moreover, if the
subgroup is separable (i.e., the intersection of the finite index
subgroups containing it), then the surface is separable \cite{scott}.
Thus our results have purely group-theoretic formulations.  In
particular, we see that many infinite 3-manifold groups have the
property that any surface subgroup is virtually abelian.  In fact, it
is easy to find examples with no Klein bottles (one must just avoid
Siefert fibered pieces with non-orientable base) and thus see that many
infinite 3-manifold groups have no non-abelian surface subgroups.  For
example, we show the fundamental group of a link of an isolated
complex surface singularity always has this property.

Some of the results of Niblo and Wise \cite{niblo-wise} and
\cite{niblo-wise-gt} are also of interest in this context.  For
example, they show that subgroup separability fails for any graph
manifold which is not a Seifert manifold or covered by a torus bundle
and they show that a non-separable horizontal surface in a graph
manifold can only be lifted to finitely many finite covers of $M^3$.

\noindent{\bf Acknowledgements}.~~This research was supported by the
Australian Research Council. I am grateful to S.K. Roushon for the
question that led to this paper and to Hyam Rubinstein for useful
conversations. 

After seeing an electronic preprint of this paper
\cite{neumann-preprint}, Buyalo brought my attention to the
interesting series of papers \cite{buyalo-kobelskii1},
\cite{buyalo-kobelskii2}, \cite{buyalo-kobelskii3},
\cite{buyalo-kobelskii4} in which he and Kobel${}'$ski\v\i{} study
existence of metrics of various types on graph manifolds. Of
particular interest is their study of what they call ``isometric
geometrizations'' in \cite{buyalo-kobelskii2}.  They cite \cite{leeb}
to say that the existence of such a geometrization is equivalent to
the existence of a metric of non-positive sectional curvature.
Their conditions are rather close to the conditions arising here and
in \cite{neumann}, but it is not clear to the author why this is so.

\section{Main results}
\label{sec:main results}

From now on we assume $M^3$ is a closed connected graph-manifold, that
is, a closed connected manifold obtained by pasting together compact
Seifert fibered $3$-manifolds along boundary components.  We are
interested in two properties:
\begin{enumerate}
\item[(I)] $M^3$ has an immersed $\pi_1$-injective surface of negative
Euler characteristic;
\item[(VE)] $M^3$ has a virtually embedded $\pi_1$-injective surface of
negative Euler characteristic (i.e., some finite cover of $M^3$ has an
embedded such surface).
\end{enumerate}
There is no loss of generality in assuming $M^3$ is irreducible, since
a $\pi_1$-injective surface can be isotoped to be disjoint from any
embedded $S^2$.  The properties (I) and (VE) are preserved on
replacing $M^3$ by a finite cover, so there is also no loss of
generality in assuming $M^3$ is orientable.  Moreover, if $M^3$ has a
cover which is a torus bundle over $S^1$ then $M^3$ fails both
properties (I) and (VE) so we may assume $M^3$ is not of this type.
It is then easy to show (see \cite{neumann} p.364) that, after passing
to a double cover if necessary, we may assume that $M^3$ can be cut
along a family of disjoint embedded $\pi_1$-injective tori into Seifert
fibered pieces $M_1,\dots, M_s$ satisfying:
\begin{itemize}
\item Each $M_i$ is Seifert fibered over
  orientable\footnote{Orientability is for convenience of proof and is
    not actually needed for our main Theorem \ref{th:main}.}
  base-orbifold of negative orbifold Euler characteristic;
\item no $M_i$ meets itself along one of the separating tori.
\end{itemize}

We may also assume that none of the separating tori is redundant.
This means:
\begin{itemize}
\item For each separating torus $T$ the fibers of the Seifert pieces
  on each side of $T$ have non-zero intersection number (which
  we denote $p(T)$) in $T$.
\end{itemize}

Each Seifert fibered piece $M_i$ has a linear foliation of its
boundary by the Seifert fibers of the adjacent Seifert fibered pieces.
The rational Euler number $e_i$ of the Seifert fibration of $M_i$ with
respect to these foliations is defined as $e_i=e(\hat M_i\to \hat
F_i)$, where $\hat M_i\to \hat F_i$ is the closed Seifert fibration
obtained by filling each boundary component of $M_i$ by a solid torus
whose meridian curves match the foliation. As in \cite{neumann} we
define the \emph{decomposition matrix} for $M^3$ to be the symmetric
matrix $A(M^3)=(A_{ij})$ with
\begin{align*}
  A_{ii}&=e_i \\
  A_{ij}&=\sum_{T\subset M_i\cap M_j}\frac{1}{|p(T)|}\quad (i\ne j),
\end{align*}
\noindent where the sum is over components $T$ of $M_i\cap M_j$ and
$p(T)$ is, as above, the intersection number in $T$ of fibers from the
two sides of $T$.

$A(M^3)$ is a symmetric rational matrix with non-negative off-diagonal
entries. Moreover, the graph on $s$ vertices, with an edge connecting
vertices $i$ and $j$ if and only if $A_{ij}\ne 0$, is a connected graph.
Given any matrix $A$ with these properties, it is easy to realise it
as $A(M^3)$ for some $M^3$.

By reordering indices we may put $A(M^3)$ in block form
\begin{equation*}
  \begin{pmatrix}
    P&Z\\
    Z^t&N
  \end{pmatrix}
\end{equation*}
where $P$ has non-negative diagonal entries and $N$ has non-positive
diagonal entries%
\footnote{Each zero diagonal entry can be put in either $P$ or $N$.
  Notation here therefore differs from \cite{neumann} where they were
collected in their own block}.  Let $P_-$ be the result of
multiplying the diagonal entries of $P$ by $-1$ and put
\begin{equation*}
  A_-(M^3):=
  \begin{pmatrix}
    P_-&Z\\ Z^t& N
  \end{pmatrix}.
\end{equation*}

\begin{theorem}\label{th:main}
$M^3$ satisfies condition (I), that
  is, $M^3$ has an immersed $\pi_1$-injective surface of negative
  Euler characteristic, if and only if either $A_-(M^3)$ has a
  positive eigenvalue or it is negative and indefinite and all
  diagonal entries of $A(M^3)$ have the same sign (in which case
  $A_-(M^3)$ is negative semidefinite and $M^3$ even satisfies (VE)).
  
  $M^3$ satisfies condition (VE), that is, $M^3$ has a virtually
  embedded $\pi_1$-injective surface of negative Euler characteristic,
  if and only if one of $P_-$ or $N$ is not negative definite.
\end{theorem}
\noindent(It is an elementary but not completely trivial exercise to see that
the algebraic condition of this theorem for (I) follows from the
algebraic condition for (VE). This also follows from the proof of the
theorem.)
\begin{example}
  If $M^3$ is the link of an isolated complex surface singularity
  then, as discussed in \cite{neumann}, $A(M^3)$ is negative definite
  (so $A(M^3)=A_-(M^3)$), so $M^3$ fails condition (I), and hence also
  fails (VE).  Since it is known by \cite{neumann-t} that the Seifert
  components of such an $M^3$ all have orientable base, $M^3$ has no
  immersed Klein bottles, so all infinite surface subgroups of
  $\pi_1(M^3)$ are abelian.
\end{example}

\begin{example}\label{ex:2}
  If $M^3=M_1\cup M_2$ has just two Seifert components, so the
  decomposition matrix is $A(M^3)=(A_{ij})_{1\le i,j\le 2}$, put
  $D:=(A_{11}A_{22})/{A_{12}^2}$.  Theorem \ref{th:main} implies:
  \begin{align*}
  M^3 \text{ satisfies (I)}\quad&\Leftrightarrow\quad-1<D\le 1;\\
  M^3 \text{ satisfies (VE)} \quad&\Leftrightarrow\quad\phantom{-}0\le
    D\le 1.
  \end{align*}
  In \cite{neumann} it was shown that this $M^3$ is virtually fibered
  over $S^1$ (i.e., has a finite cover that is fibered) if and only if
  either $0<D\le1$ or $A_{11}=A_{22}=0$.  Moreover, $M^3$ itself
  fibers over $S^1$ if and only if $D=1$.  The manifolds of this
  example were classified up to commensurability in \cite{neumann} by
  two rational invariants, one of which is the above $D$.
\end{example}

One can ask also about compact graph manifolds $M^3$ with non-empty
boundary.  If we assume $M^3$ is orientable, irreducible and not one
of the trivial cases $D^2\times S^1, T^2\times I$, or $I$-bundle over
the Klein bottle then $M^3$ always has virtually embedded surfaces of
negative Euler characteristic by \cite {rubinstein-wang}.  In fact,
Wang and Yu \cite{wang-yu} show more: $M^3$ is virtually fibered over
$S^1$.  This can also be deduced using only matrix algebra (but a
little effort) from \cite{neumann}, where a necessary and sufficient
condition for virtual fibering of a closed graph manifold is given in
terms of the decomposition matrix $A(M^3)$. This approach actually
proves the stronger result (we omit details):
\begin{theorem}
  If $M^3$ is an oriented irreducible graph manifold with nonempty
  boundary then there exists a fibration $\partial M^3\to S^1$ which
  extends to a virtual fibration of $M^3$ to $S^1$ (that is, each
  fiber of the virtual fibration is a virtually embedded surface whose
  boundary is parallel to the given fibration of $\partial M^3$).\qed
\end{theorem}

\section{Proofs}
\label{sec:proofs}

The necessary and sufficient condition for condition (VE) was proved
in \cite{neumann}, so in this paper we just prove the condition for
(I).  We start with a discussion of Seifert fibered manifolds.

If $\pi\colon M\to F$ is a Seifert fibration, a proper immersion
$f\colon S\to M$ of a surface $S$ is \emph{horizontal} if it is
transverse to all fibers of $\pi$. Equivalently, $\pi\circ f$ is a
covering map of $S$ to the orbifold $F$.

Suppose $\pi\colon M\to F$ is a Seifert fibration with $F$ connected
and orientable and $\chiorb(F)<0$ (orbifold Euler characteristic).
Assume $M$ has non-empty boundary.
On each torus $T\subset\partial M$ let a section $m_T$ to the
Seifert fibration be given.  Then $e(M\to F)$ is defined with respect
to these sections. We orient each $m_T$ consistently with $\partial
F$.  Let $f\colon S\to M$ be a horizontal immersion of a surface $S$.
Orient $S$ so $\pi\circ f$ preserves orientation.  Denote the boundary
components of $S$ that lie in $T$ by $c_{T1},\dots,c_{T k_T}$.
Using $f$ to denote a generic fiber of $\pi$ we have integers $a_{T
  \beta }, b_{T \beta }$, $\beta=1,\dots,k_T$, so that the
following homology relations hold:
\begin{equation}
  \label{eq:homol}
  [c_{T \beta }]=a_{T \beta }[m_{T }]+b_{T \beta
  }[f]\in H_1(T ),\quad a_{T\beta}>0.
\end{equation}
\begin{lemma}\label{le:seifert}
  If $a$ is the degree of $\pi\circ f\colon S\to F$ then
  \begin{equation}\label{eq:l1}
    \sum_{\beta =1}^{k_T }a_{T \beta }=a\quad
    \text{for each $T\subset\partial M$.}
  \end{equation}
  Moreover,
  \begin{equation}\label{eq:l2}
    \sum_{T \subset\partial M}\sum_{\beta =1}^{k_T }b_{T \beta }=ae.
  \end{equation}
  Conversely, suppose that for each boundary component $T_\beta $ 
  there is given a family $c_{T 1},\dots,c_{T k_T }$ of immersed 
  curves transverse to the fibers of $\pi$ satisfying homology
  relations (\ref{eq:homol}), and that equations (\ref{eq:l1}) and
  (\ref{eq:l2}) are satisfied.  Then there exist integers $d_0>0,
  n_0>0$ so that for any positive integer multiple $d$ of $d_0$ and
  $n$ of $n_0$ the family of curves $c_{T \beta \gamma}^n$,
  $T\subset\partial M$, $\beta =1,\dots,k_T $, $\gamma=1,\dots,d$,
  obtained as follows, bounds an immersed horizontal surface.  For
  $\gamma=1,\dots,d$ we take $c_{T \beta \gamma}^n$ as a copy of
  the immersed curve obtained by going $n$ times around the curve
  $c_{T \beta }$.
\end{lemma}
\begin{proof}
  The left side of equation (\ref{eq:l1}) is $\degree(\pi\circ
  f|(S\cap T )\colon S\cap T \to \partial_T F)$,
  which is the degree of $\pi\circ f$, proving (\ref{eq:l1}). Now the
  sum over all $T$ and $\beta $ of the curves $c_{T \beta }$ is
  null-homologous in $M$ and using equation (\ref{eq:l1}) this says
  $\sum_T a[m_T ]+ \bigl(\sum_{T ,\beta }b_{T
    \beta }\bigr)[f]=0$.  Equation (\ref{eq:l2}) now follows from the
  fact that $e$ can be defined by the equation $\sum_{T
    =1}^l[m_T ]=-e[f]$ in $H_1(M;\Q)$.
  
  For the converse we first apply Lemma 5.1 of \cite{neumann} which
  implies (taking $d_0e$ to be integral in that lemma):
  
  There are positive integers $d_0$ and $n_0$ such that $d_0e\in \Z$
  and such that for any multiples $d$ of $d_0$ and $n$ of $n_0$ there
  is a covering $p\colon M'\to M$ satisfying
  \begin{itemize}
  \item The lifted Seifert fibration of $M'$ has no singular fibers
    (so $M'\cong F'\times S^1$, where $F'$ is the base surface for
    the fibration of $M'$),
  \item $p$ has degree $dn^2$,
  \item each boundary torus $T $ of $M$ is covered by $d$
    boundary tori $T_{\gamma}$, $\gamma=1,\dots,d$ of $M'$,
    each of which is a copy of the unique connected
    $(\Z/n\times\Z/n)$-cover of $T $.
  \end{itemize}
  Now each curve $c_{T \beta }$ lifts to $n$ curves in $T_{\gamma}$,
  each still of slope $a_{T \beta }/b_{T \beta }$.  Pick one of these
  and call it $c_{T \beta \gamma}$. If we can find a horizontal
  surface $S'$ in $M'$ spanning the family of curves $\{c_{T \beta
    \gamma}:T\subset\partial M, \beta =1,\dots,k_T ,
  \gamma=1,\dots,d\}$, then its image in $M$ is the desired surface.
  
  The identification of $M'$ with $F'\times S^1$ gives 
  meridian curves $m'_{T \gamma}\in T_{\gamma}$ and with
  respect to these the Euler number of $M'\to F'$ is $0$.  Thus the
  curves $c_{T \beta \gamma}$ satisfy homology relations
  $c_{T \beta \gamma}=a_{T \beta }[m'_{T
    \gamma}]+b'_{T \beta }[f']$ for some $b'_{T \beta }$
  with $\sum_{T ,\beta }b'_{T \beta }=0$.  We are thus
  looking for a connected surface $S'$ mapping to $M'\times S^1$ by a
  map $(g,h)\colon S'\to F'\times S^1$ such that:
\begin{itemize}
\item the map $g$ is a covering of degree $da$ and the boundary
  component corresponding to $T_\gamma$ of $F'$ is covered by exactly
  $k_T $ boundary components $\partial_{T \beta \gamma}S', \beta
  =1,\dots,k_T $ of $S'$, with degrees $a_{T 1},\dots,a_{T k_T }$;
\item the map $h$ has degree $b'_{T \beta }$ on $\partial_{T
    \beta \gamma}S'$.
\end{itemize}
If $S'$ is
connected then the fact that $[S',S^1]=H^1(S';\Z)$ and the exact cohomology
sequence
\begin{equation*}
  H^1(S';\Z)\to H^1(\partial S';\Z)\to H^2(S',\partial S';\Z)=\Z.
\end{equation*}
shows that $h\colon S'\to S^1$ exists with degree $b'_{T\beta}$
on each $\partial_{T\beta\gamma}S'$ if and only if $\sum
b'_{T\beta}=0$.  Thus the only issue is finding a connected cover
$S'$ of $F'$ with $g\colon S'\to F'$ as above.

Since $F'$ is a $dn$-fold cover of the orbifold $F$ and $F'$ has $ld$
boundary components, where $l$ is the number of boundary components of
$M$, we have $2-2\genus(F')=dn\chiorb(F)+dl$, so $\genus(F')>0$ as
soon as $n$ is chosen large enough.  We therefore assume
$\genus(F')>0$.  We also choose $d_0$ even. The existence of a
connected cover with prescribed degree and boundary behaviour then
follows from the following lemma, since the parity condition of the
lemma is $da\chi(F')\equiv d\sum_T k_T $ (mod $2$).
\end{proof}
\begin{lemma}
  If $F'$ is an orientable surface of positive genus and a degree
  $\alpha\ge1$ is specified and for each boundary component a
  collection of degrees summing to $\alpha$ is also specified, then a
  connected $\alpha$-fold covering $S'$ of $F'$ exists with prescribed
  degrees on the boundary components over each boundary component of
  $F'$ if and only if the prescribed number of boundary components of
  the cover has the same parity as $\alpha\chi(F')$.
\end{lemma}
\begin{proof}
  This lemma appears to be well known, although weaker results have
  appeared several times in the literature.  It is assumed implicitly
  in the proof of Lemma 2.2 of \cite{rubinstein-wang} (which has a
  minor error, since the parity condition is overlooked).  The parity
  condition arises because the Euler characteristic of a compact
  orientable surface has the same parity as the number of its boundary
  components.  Alternatively, if one represents the cover by a
  homomorphism of $\pi_1(F')$ to the symmetric group $Sym_\alpha$ of a
  fiber, the parity condition arises because the product of the
  permutations represented by boundary components is a product of
  commutators and is hence an even permutation.  The existence of
  $S'\to F'$ with the given constraints can be seen by constructing a
  homomorphism of $\pi_1(F')$ to $Sym_\alpha$ with transitive image
  which maps the boundary curves to permutations with the desired
  cycle structure.  Such a homomorphism exists by the result of
  Jacques et al.\ \cite{jacques-et-al} that any even permutation on
  $n$ symbols is a commutator of an $n$-cycle and an involution.
\end{proof}

We now return to the graph manifold $M^3$ of Section 
\ref{sec:main results} 
which is glued together from Seifert fibered manifolds
$M_1, \dots M_s$.  For each $M_i$ we choose an orientation of the base
surface of the Seifert fibration.  We can assume we have done this so
that for each separating torus $T$ the intersection number $p(T)$ of
the Seifert fibers from the two sides of $T$ is positive.  Indeed, if
this is not possible, then, as pointed out in \cite{neumann}, we can
replace $M^3$ by a commensurable graph manifold $M'$ with the same
decomposition matrix for which it is possible (in fact $M^3$ and $M'$
have a common 2-fold cover).  From now on we will therefore assume all
$p(T)$ are positive.

We first prove that $A_-(M^3)$ not being negative definite is
necessary for having an immersed $\pi_1$-injective surface $S$ in
$M^3$ of negative Euler characteristic.  As is proven in
\cite{rubinstein-wang} (Lemma 3.3), such a surface is homotopic to an
immersed surface whose intersection with each $M_i$ consists of a
union of horizontal surfaces and possibly also some $\pi_1$-injective
vertical annuli.  We will therefore assume that our surface has
already been put in this position.  For the moment we assume also, for
simplicity, that $S$ is horizontal in $M^3$, that is, vertical annuli
do not occur.

Fix an index $i$ and consider the intersection of our immersed surface
$S$ with $M_i$.  We orient this immersed surface in $M_i$ so that it
maps orientation preservingly to the base surface of $M_i$.  We also
choose meridian curves in the boundary tori of $M_i$ and thus obtain a
collection of integer pairs $(a_{T\beta}, b_{T\beta})$ as in
Lemma \ref{le:seifert} satisfying the relations of that lemma.  Note
that the $e$ in that lemma is not $e_i$, since it is Euler number with
respect to the chosen meridians rather than with respect to the
Seifert fibers of neighbouring Seifert fibered pieces to $M_i$. We
denote it therefore $e'_i$. We denote the degree $a$ appearing in the
lemma by $a_i$.

Our orientation of $S\cap M_i$ induces an orientation on each boundary
curve of this surface.  Each such curve also inherits an orientation
from the piece of surface it bounds in a neighbouring
Seifert fibered piece. Call a curve \emph{consistent} if these two
orientations agree. For fixed $T$ denote by $a^+_T$ the sum
of the $a_{T\beta}$'s corresponding to consistent curves and
$a^-_T$ the sum of the remaining $a_{T\beta}$'s.  Define
$b^+_T$ and $b^-_T$ similarly.  Thus equations (\ref{eq:l1})
and (\ref{eq:l2}) become
\begin{align}
  a^+_T+a^-_T&=a_i\quad\text{for each $T\subset\partial M_i$},
  \label{eq:1}\\
  \sum_{T\subset\partial M_i}(b^+_T+b^-_T)&=a_ie'_i.\label{eq:2}
\end{align}

For given $T\subset\partial M_i$ we denote by $T'$ the same torus
considered as a boundary component of the Seifert piece $M_j$ adjacent
to $M_i$ across $T$. The pair $(a^+_T,b^+_T)$ gives coordinates of the
homology class represented by the sum of the consistent curves in $T$
with respect to the basis of $H_1(T)$ coming from meridian and fiber
in $M_i$.  The same homology class will be given by a pair
$(a^+_{T'},b^+_{T'})$ with respect to meridian and fiber in $M_j$ with
\begin{equation}\label{eq:m1}
  \begin{pmatrix}
    a^+_{T'}\\b^+_{T'}
  \end{pmatrix}
=
\begin{pmatrix}
  q(T)&p(T)\\-p'(T)&-q'(T)
\end{pmatrix}
\begin{pmatrix}
  a^+_T\\b^+_T
\end{pmatrix},
\end{equation}
where the square matrix is the appropriate change-of-basis matrix.
Our notation for this matrix agrees with page 366 of \cite{neumann};
in particular, $p(T)$ has its meaning of intersection number of fibers
of $M_i$ and $M_j$ in $T$. The matrix has determinant $-1$, since $T$
has opposite orientations viewed from $M_i$ and $M_j$.  We also have:
\begin{equation}\label{eq:m2}
  \begin{pmatrix}
    a^-_{T'}\\b^-_{T'}
  \end{pmatrix}
= -
\begin{pmatrix}
  q(T)&p(T)\\-p'(T)&-q'(T)
\end{pmatrix}
\begin{pmatrix}
  a^-_T\\b^-_T
\end{pmatrix}.
\end{equation}
The first entries of matrix equations (\ref{eq:m1}) and (\ref{eq:m2})
are the equations 
\begin{equation}\label{eq:first}
  a^\pm_{T'}=\pm\bigl(q(T)a^\pm_T+p(T)b^\pm_T\bigr)
\end{equation}
that we
can solve for $b^\pm_T$ in terms of $a^\pm_T$ and
$a^\pm_{T'}$ to give:
\begin{equation}\label{eq:b}
b^\pm_T=(\pm a^\pm_{T'}-q(T)a^\pm_T)/p(T).
\end{equation}
Equation (\ref{eq:2}) thus becomes:
\begin{equation}\label{eq:2o}
  \sum_{T\subset\partial
  M_i}\left(\frac{a^+_{T'}-q(T)a^+_T}{p(T)}+\frac{-a^-_{T'}-q(T)a^-_T}{p(T)} 
  \right)=a_ie'_i.
\end{equation}
Using equation (\ref{eq:1}) this becomes
\begin{equation}\label{eq:2a}
  \sum_{T\subset\partial
  M_i}\frac{a^+_{T'}-a^-_{T'}}{p(T)}=a_i\biggl(e'_i+\sum_{T\subset\partial
  M_i}\frac{q(T)}{p(T)}\biggr).
\end{equation}
As discussed on page 366 of \cite{neumann}, $q(T)/p(T)$ is the change
of Euler number $e(M_i\to F_i)$ on replacing the meridian at $T$ by
the fiber of $M_j$.  Thus the right side of (\ref{eq:2a}) is $a_ie_i$,
so equation (\ref{eq:2a}) says
\begin{equation}
  \label{eq:2b}
  \sum_{T\subset\partial
    M_i}\frac{a^+_{T'}-a^-_{T'}}{p(T)}=a_ie_i.
\end{equation}

Consider the summands on the left with $T\subset\partial
M_i\cap\partial M_j$.  Since $a^+_{T'}+a^-_{T'}=a_j$ and $a^+_{T'}$
and $a^-_{T'}$ are both non-negative, each summand is no larger in
magnitude than $a_j/p(T)$. Their sum is therefore no larger in
magnitude than 
\begin{equation*}
  a_j\left(\sum_{T\subset\partial M_i\cap\partial
  M_j}\frac1{|p(T)|}\right)=a_jA_{ij}.
\end{equation*}
  We write their sum therefore as
$-a_jA'_{ij}$ with $|A'_{ij}|\le A_{ij}$, so (\ref{eq:2b}) becomes
\begin{equation}
  \label{eq:2c}
  -\sum_{j\ne i}A'_{ij}a_j=a_ie_i.
\end{equation}
Recalling that $e_i=A_{ii}$ and putting $A'_{ii}=A_{ii}$ we can
rewrite this as
\begin{equation}
  \label{eq:2d}
  \sum_{j=1}^sA'_{ij}a_j=0.
\end{equation}
We have thus shown that the decomposition matrix $A(M^3)$ has the
property that it can be made to have non-trivial kernel by replacing
each off-diagonal entry by some rational number of no larger
magnitude.  The fact that $A_-(M^3)$ is not negative definite thus
follows from the following lemma.

If $A$ is a matrix with non-negative off-diagonal entries
then we will use the term \emph{reduction of $A$} for a matrix $A'$
with $|A'_{ij}|\le A_{ij}$ for all $i\ne j$ and $A'_{ii}=A_{ii}$ for
all $i$. 
\begin{lemma}\label{le:matrix}
  Let $A=(A_{ij})$ be a square symmetric matrix over $\Q$ with
  $A_{ij}\ge0$ for $i\ne j$.  Then there exists a (not necessarily
  symmetric) singular rational reduction $A'=(A'_{ij})$ of $A$ if and
  only if the matrix $A_-$ (obtained by replacing each positive
  diagonal entry of $A$ by its negative) is not negative definite.
  Moreover such an $A'$ can then be found which annihilates a non-zero
  vector with non-negative entries.
\end{lemma}
We postpone the proof of this Lemma and first return to the proof of
Theorem \ref{th:main}.  The fact that $A_-(M^3)$ is not negative
definite is not quite proved, since we assumed vertical annuli do not
exist in our $\pi_1$-injective surface. If we do have vertical annuli
we choose orientations on them. Then we can characterise their
boundary components as consistent or non-consistent as before.
Equations (\ref{eq:1}) and (\ref{eq:2}) then still hold, so the above
proof goes through unchanged.

For the converse, suppose that the decomposition matrix
$A(M^3)=(A_{ij})$ is not negative.  We shall show that this actually
implies the existence of a horizontal surface (i.e., with no vertical
annuli).  Our condition on $A_-(M^3)$ is that it has a positive
eigenvalue, which is an open condition, so we can reduce each non-zero
off-diagonal entry slightly without changing it.  By the above lemma
we can thus assume there exists a rational matrix $(A'_{ij})$ with
$|A'_{ij}|< A_{ij}$ for each $i\ne j$ with $A_{ij}\ne 0$ and with
$A'_{ii}=A_{ii}$ for each $i$ such that equation (\ref{eq:2d}) (or the
equivalent equation (\ref{eq:2c})) holds for some non-zero vector
$(a_1,\dots,a_s)$ with non-negative rational entries.  For each $i\ne
j$ we then define $a^+_{T'}$ and $a^-_{T'}$, for each boundary torus
$T'$ of $M_j$ that lies in $M_i\cap M_j$, by the equations
\begin{align*}
  a^+_{T'}&=\frac{A_{ij}-A'_{ij}}{2A_{ij}}\,a_j\\
  a^-_{T'}&=\frac{A_{ij}+A'_{ij}}{2A_{ij}}\,a_j.
\end{align*}
Note that these imply that $a^\pm_{T'}>0$ whenever $a_j\ne0$ and
\begin{align*}
  a^+_{T'}+a^-_{T'}&=a_j\\
  a^+_{T'}-a^-_{T'}&=-(A'_{ij}/A_{ij})a_j.
\end{align*}
Thus, equation (\ref{eq:1}) holds, and, working backwards via
equations (\ref{eq:2b}), (\ref{eq:2a}) and (\ref{eq:2o}) we see that
(\ref{eq:2}) holds if we define $b^\pm_{T}$ by equation (\ref{eq:b}).
Moreover, by multiplying our original vector $(a_j)$ by a suitable
positive integer we may assume that the $a^\pm_T$ and $b^\pm_T$ are
all integral.

Now, (\ref{eq:b}) is equivalent to (\ref{eq:first}) which can also be
written
\begin{equation}\label{eq:next1}
  a^\pm_T=\pm\bigl(q(T')a^\pm_{T'}+p(T')b^\pm_{T'}\bigr),
\end{equation}
by exchanging the roles of $T$ and $T'$. But $q(T')=q'(T)$
and $p(T')=p(T)$. In fact
\begin{equation}
  \begin{pmatrix}
      q(T')&p(T')\\-p'(T')&-q'(T')
  \end{pmatrix}
=
  \begin{pmatrix}
      q'(T)&p(T)\\-p'(T)&-q(T)
  \end{pmatrix},
\end{equation}
since the coordinate change matrix for $T'$ is the inverse of the one
for $T$.  Thus (\ref{eq:next1}) implies
\begin{equation}\label{eq:next}
  a^\pm_T=\pm\bigl(q'(T)a^\pm_{T'}+p(T)b^\pm_{T'}\bigr).
\end{equation}
Inserting (\ref{eq:first}) in
(\ref{eq:next}) and simplifying, using the fact that
$1-q'(T)q(T)=-p'(T)p(T)$, gives
$p(T)b^\pm_{T'}=\pm\bigl(p'(T)p(T)a^\pm_T+q'(T)p(T)b^\pm_T\bigr)$,
whence
\begin{equation}
  b^\pm_{T'}=\pm\bigl(p'(T)a^\pm_T+q'(T)b^\pm_T\bigr).
\end{equation}
With equation (\ref{eq:first}) this gives the matrix equations
(\ref{eq:m1}) and (\ref{eq:m2}) which imply that the curve $c^\pm_T$
in $T$ defined by coordinates $(a^\pm_T,b^\pm_T)$ with respect to
meridian and fiber in $M_i$ is the same as the curve in $T'$ defined
by $(a^\pm_{T'},b^\pm_{T'})$ with respect to meridian and fiber in
$M_j$. We thus have a pair of curves in each separating torus so that
the curves in the boundary tori of each Seifert piece $M_i$ satisfy
the numerical conditions of Lemma \ref{le:seifert}.  By that Lemma, we
can find $d$ and $n$ so that if we replace each of the curves $c$ in
question by $d$ copies of the curve $c^n$, then the curves span a
horizontal surface in each $M_i$. These surfaces fit together to give
the desired surface in $M^3$.

It remains to discuss the case that $A_-(M^3)$ is negative
but not definite.  We postpone this until after the proof of the lemma.

\begin{proof}[Proof of Lemma \ref{le:matrix}]
  We first note that if $A$ has a singular reduction then it has a
  reduction that annihilates a vector with non-negative entries.
  Indeed, if we have a reduction $A'$ that annihilates the non-trivial
  vector $(x_i)$, then for each $i$ with $x_i<0$ we multiply the
  $i$-th row and column of $A'$ by $-1$.  The result is a reduction
  $A''$ which annihilates $(|x_i|)$.  We next note that the property
  of $A$ having a singular reduction is unchanged if we change the
  sign of any diagonal entry of $A$, since if $A'$ is a singular
  reduction for $A$ then multiplying the corresponding row of $A'$ by
  $-1$ gives a singular reduction of the modified matrix.  Thus, we
  may assume without loss of generality that our initial matrix $A$
  has non-positive diagonal entries.
  
  Suppose $A$ is symmetric with non-negative off-diagonal entries and
  non-positive diagonal entries and suppose $A$ has a singular
  reduction $A'$, say $A'x=0$ with $x$ a non-zero vector.  Then
  $x^t(A'+(A')^t)x=0$, so $\frac12(A'+(A')^t)$ is an indefinite
  symmetric reduction of $A$.  In \cite{neumann} it is shown that a
  symmetric reduction of a negative definite matrix with non-negative
  off-diagonal entries is again negative definite. Thus $A$ is not
  negative definite.
  
  Conversely, suppose $A$ is a rational symmetric matrix with
  non-negative off-diagonal entries and non-positive diagonal entries
  and suppose $A$ is not negative definite. We want to show the
  existence of a singular rational reduction of $A$. If $A$ is itself
  singular we are done, so we assume $A$ is non-singular. Assume first
  that only one eigenvalue of $A$ is positive. Consider a piecewise
  linear path in the space of reductions of $A$ that starts with $A$
  and reduces each off-diagonal entry to zero, one after another. This
  path ends with the diagonal matrix obtained by making all
  off-diagonal entries zero, which has only negative eigenvalues, so
  the determinant of $A$ changes sign along this path. It is thus zero
  at some point of the path.  Since determinant is a linear function
  of each entry of the matrix, the first point where determinant is
  zero is at a matrix with rational entries.  We have thus found a
  rational singular reduction of $A$.  If $A$ has more than one
  positive eigenvalue, consider the smallest principal minor of $A$
  with just one non-negative eigenvalue. First reduce all off-diagonal
  entries that are not in this minor to zero and then apply the above
  argument just to this minor.
\end{proof}

To complete the proof of Theorem \ref{th:main} we must discuss the
case that $A_-(M^3)$ is negative indefinite. We need some algebraic
preparation.

Let $A$ be a symmetric $s\times s$ matrix. The $s$-vertex graph with
an edge connecting vertices $i$ and $j$ if and only if $A_{ij}\ne 0$
will be called the \emph{graph of $A$}. The submatrices of $A$
corresponding to components of this graph will be called the
\emph{components} of $A$.  By reordering rows and columns, $A$ can be
put in block diagonal form with its components as the diagonal blocks.
If $A$ has just one component we call $A$ \emph{connected}.

\begin{proposition}
  Let $A$ be a symmetric $s\times s$ matrix with non-negative
  off-diagonal entries such that $A$ is connected. Then $A$ is
  negative if and only if there exists a vector $\bfa=(a_j)$ with
  positive entries such that $A\bfa$ has non-positive entries.
  Moreover, in this case $A$ is negative definite unless $A\bfa=0$, in
  which case $\bfa$ generates the kernel of $A$.
\end{proposition}
\begin{proof}
  Suppose $\bfa$ has positive entries. For any vector $\bfx=(x_j)$ we
  can write
  \begin{equation}
    \label{eq:bilin}
    \bfx^tA\bfx=\sum_ia_i\left(\sum_j A_{ij}a_j\right)
\left(\frac{x_i}{a_i}\right)^2+
\sum_{i<j}\bigl(-A_{ij}a_ia_j\bigr)\left(\frac{x_i}{a_i}-\frac{x_j}{a_j}\right)^2.
  \end{equation}
  If $A\bfa$ has non-positive entries then both terms on the right are
  clearly non-positive, proving that $A$ is negative.  Moreover, since
  the $s$-vertex graph determined by nonvanishing of $A_{ij}$ is
  connected, the second term on the right vanishes if and only if
  $x_i/a_i=x_j/a_j$ for all $i,j$, that is, $\bfx$ is a multiple of
  $\bfa$. In this case the first term vanishes if and only if $\bfx=0$
  or $A\bfa=0$.
  
  Conversely, suppose $A$ is a symmetric negative matrix with
  non-negative off-diagonal entries.  Then its diagonal entries are
  non-positive, and if any diagonal entry is zero then all other
  entries in the corresponding row and column must be zero. If a
  diagonal entry is non-zero, then, since it is negative, we can add
  positive multiples of the row and column containing it to other rows
  and columns, to make zero all off-diagonal entries in its row and
  column.  This preserves the properties of $A$ of being a symmetric
  negative matrix with non-negative off-diagonal entries.  It thus
  follows that we can reduce $A$ to a diagonal matrix using only
  ``positive'' simultaneous row and column operations, so we have
  $P^tAP=D$ where $P$ is invertible with only non-negative entries and
  $D$ is diagonal.  If $A$ is non-singular then $A^{-1}$=$PD^{-1}P^t$
  and this is a matrix with non-positive entries. Thus the negative
  sum of the columns of $A^{-1}$ is a vector $\bfa$ with positive
  entries and $A\bfa=(-1,\dots,-1)^t$, so $\bfa$ is as required.  If
  $A$ is singular, then $D$ has a zero entry, and the corresponding
  column of $P$ is a non-trivial vector $\bfa$ with non-negative
  entries such that $A\bfa=0$. Thus, in this case $\bfa$ is as
  required if we show that it has no zero entries.  Suppose $\bfa$ did
  have zero entries.  By permuting rows and columns of $A$ we can
  assume they are the last few entries of $\bfa$, so
  $\bfa=\begin{pmatrix}\bfa_0\\{\bf0}
  \end{pmatrix}$ with
  no zero entries in $\bfa_0$, and $A$ has block form
  $\begin{pmatrix} A_0&B^t\\B&A_1
  \end{pmatrix}$ with $\begin{pmatrix} A_0\\B
  \end{pmatrix}\bfa_0=0$.  Since 
  $B$ has non-negative entries and $\bfa_0$ has only positive entries,
  this implies $B=0$. This contradicts the connectedness of $A$.
\end{proof}
\begin{corollary}\label{co:negdef}
  Suppose $A$ is a symmetric connected negative matrix with
  non-negative off-diagonal entries. If $A^s$ is a symmetric reduction
  of $A$ such that some off-diagonal entry has been reduced in
  absolute value then $A^s$ is negative definite.
\end{corollary}
\begin{proof}
  Let
  $\bfa$ be a vector with positive entries such that $A\bfa$ has
  non-positive entries.  
  
  Suppose first that $A^s$ has non-negative off-diagonal entries.  Then
  $A^s\bfa$ has non-positive entries and is non-zero, so $A^s$ is
  negative definite by the preceding proposition.  In general, let
  $A'$ be the reduction of $A$ with $A'_{ij}=|A^s_{ij}|$ for $i\ne
  j$.  Then $A'$ is negative definite by what has just been said, and
  $A^s$ is a reduction of $A'$, so $A^s$ is negative definite by
  \cite{neumann}.
\end{proof}
\begin{corollary}\label{co:home}
  Suppose $A$ is a symmetric connected matrix with non-negative
  off-diagonal entries for which $A_-$ (the result of multiplying
  positive diagonal entries of $A$ by $-1$) is negative indefinite.
  Let $A'$ be a singular reduction of $A$ and $A''$ the result of
  multiplying each row of $A'$ with positive diagonal entry by $-1$.
  Then $A''$ is symmetric and it satisfies $|A''_{ij}|=A_{ij}$ for all
  $i\ne j$.
\end{corollary}
\begin{proof}
   Since $A''$ is singular, its symmetrization
$A^s=\frac12(A''+(A'')^t)$ is not negative definite (since $A''\bfx=0$
implies $\bfx^tA^s\bfx=0$).  By Corollary \ref{co:negdef} the entries
of $A^s$ are therefore the same in absolute value as the entries of
$A$.  This implies that $A''$ was already symmetric and its entries
are the same in absolute value as those of $A$.
\end{proof}

Suppose now that $A_-=A_-(M^3)$ is negative indefinite and $M^3$
satisfies condition (I).  We want to show that the block decomposition
\begin{equation*}
A=  \begin{pmatrix}
    P&Z\\
    Z^t&N
  \end{pmatrix}
\end{equation*}
of $A$ is trivial. Suppose we have a reduction $A'$ of $A=A(M^3)$ is
realised by a $\pi_1$-injective surface as in the proof of the
necessary condition of the main theorem.  Let $(a_j)$ be as in that
proof, so it is a non-trivial vector with non-negative integer entries
which $A'$ annihilates.  Let $A_{ij}$ be a non-zero entry of the block
$Z$.  Then Corollary \ref{co:home} implies that $A'_{ij}=-A'_{ji}=\pm
A_{ij}$.  If $a_j=0$ we could replace $A'_{ij}$ by zero, which is
impossible by Corollary \ref{co:home}, so we may assume $a_j>0$.  The
condition $|A'_{ij}|=A_{ij}$ implies that either all the $a^+_T$ with
$T\subset M_i\cap M_j$ are zero (if $A'_{ij}=A_{ij}$) or all the
$a^-_T$ with $T\subset M_i\cap M_j$ are zero (if $A'_{ij}=-A_{ij}$).
The fact that $A'_{ji}= -A'_{ij}$ implies that the corresponding
$b^\pm_T$'s are not zero.  Such $(a^\pm_T,b^\pm_T)=(0,b^\pm_T)$ must
come from vertical annuli.  The fiber coordinates of boundary
components of vertical annuli sum to zero.  Using equation
(\ref{eq:b}) this says
\begin{equation*}
\sum A'_{ij}a_j=0,\quad\text{sum over $j$ with $A_{ij}$ in $Z$.}
\end{equation*}
Subtracting this equation from equation (\ref{eq:2d}) we see that the
reduction of $A$ obtained by replacing the $A'_{ij}$ corresponding to
entries of $Z$ by zero also annihilates the vector $(a_j)$.  This
contradicts Corollary \ref{co:home}, so the block decomposition of $A$
was trivial.

Conversely, if the above block decomposition of $A$ is trivial, that
is, either $A=N$ or $A=P$, then $M^3$ satisfies (VE), so it certainly
satisfies (I).\qed

\Addressesr

\end{document}

%% file: agtout.tex
\def\ifplaintex{\expandafter\ifx\csname documentclass\endcsname\relax}

\def\gtp{{\mathsurround=0pt\it $\cal G\mskip-2mu$eometry \&\ 
$\cal T\!\!$opology $\cal P\!$ublications}}  

\def\Addressesr{\bigskip
{\small \parskip 0pt \leftskip 0pt \rightskip 0pt plus 1fil \def\\{\par}
\sl\theaddress\par
\medskip
\rm Email:\stdspace\tt\theemail\hfill\rm Received:\qua\receiveddate \par}}

\def\recd{{\small Received:\qua\receiveddate\ifx\reviseddate\relax
\else\qquad Revised:\qua\reviseddate\fi\par}} 

\def\lognumber#1{\def\thelognumber{#1}}
\def\volumenumber#1{\def\thevolumenumber{#1}}
\def\volumeyear#1{\def\thevolumeyear{#1}}
\def\papernumber#1{\def\thepapernumber{#1}}
\def\pagenumbers#1#2{\def\startpage{#1}\def\finishpage{#2}}
\def\published#1{\def\publishdate{#1}}

\def\received#1{\def\receiveddate{#1}}

\def\accepted#1{\def\accepteddate{#1}}
\def\asciititle#1{\def\theasciititle{#1}}

\long\def\asciiabstract#1{\long\def\theasciiabstract{#1}}
\def\asciikeywords#1{\def\theasciikeywords{#1}}

\let\\\par\let\thelognumber\relax\let\thevolumenumber\relax
\let\thepapernumber\relax\let\thevolumeyear\relax\let\startpage\relax
\let\finishpage\relax\let\publishdate\relax\let\receiveddate\relax
\let\reviseddate\relax\let\accepteddate\relax\let\theasciititle\relax
\let\theasciiauthors\relax
\let\theasciiabstract\relax\let\theasciikeywords\relax

\let\theasciiemail\relax

\ifplaintex
\font\logobig=cmssbx10 scaled 3836
\font\logomed=cmssbx10 scaled 2557
\else
\font\logobig=cmssbx10 scaled 4200
\font\logomed=cmssbx10 scaled 2800
\fi

\long\def\makeagttitle{   
\count0=\startpage
\agt\hfill      
\hbox to 45truept{\vbox to 0pt{\vglue -13truept{\logomed A\kern -.37em{\logobig 
T}\kern -.38em G}\vss}\hss}
\break
{\small Volume \thevolumenumber\ (\thevolumeyear)
\startpage--\finishpage\nl
Published: \publishdate}

\vglue .25truein

{\parskip=0pt\leftskip 0pt plus
1fil\def\\{\par\smallskip}{\Large\bf\thetitle}\par\medskip} \vglue
0.05truein

{\parskip=0pt\leftskip 0pt plus 1fil\def\\{\par}{\sc\theauthors}
\par\medskip}%
 
\vglue 0.03truein 

{\small\leftskip 25truept\rightskip 25truept{\bf Abstract}\stdspace\theabstract

{\bf AMS Classification}\stdspace\theprimaryclass
\ifx\thesecondaryclass\relax\else; \thesecondaryclass\fi\par
{\bf Keywords}\stdspace \thekeywords\par}\vglue 7truept

}   

\ifplaintex
\font\phead=cmsl9 scaled 950
\font\pnum=cmbx10 scaled 913
\font\pfoot=cmsl9 scaled 950
\headline{\vbox to 0pt{\vskip -4.5mm\line{\small\phead\ifnum
\count0=\startpage ISSN 1472-2739 (on-line) 1472-2747 (printed)
\hfill {\pnum\folio}\else\ifodd\count0\def\\{ }%
\ifx\theshorttitle\relax\thetitle\else\theshorttitle\fi\hfill{\pnum\folio}
\else\def\\{ and }{\pnum\folio}\hfill\ifx\theshortauthors\relax\theauthors
\else\theshortauthors\fi\fi\fi}\vss}}
\footline{\vbox to 0pt{\vglue 0mm\line{\small\pfoot\ifnum\count0=\startpage
\copyright\ \gtp\hfill\else
\agt, Volume \thevolumenumber\ (\thevolumeyear)\hfill\fi}\vss}}
\else
\font=cmsl9 scaled 1050
\font\lnum=cmbx10 
\font=cmsl9 scaled 1050
\makeatletter
\def\@oddhead{{\small\lhead\ifnum\count0=\startpage ISSN 1472-2739 
(on-line) 1472-2747 (printed)\hfill {\lnum\number\count0}\else\ifodd\count0
\def\\{ }\ifx\theshorttitle\relax \thetitle \else\theshorttitle\fi\hfill
{\lnum\number\count0}\else\def\\{ and }{\lnum\number\count0}
\hfill\ifx\theshortauthors\relax 
\theauthors\else\theshortauthors\fi\fi\fi}}\def\@evenhead{\@oddhead}
\def\@oddfoot{\small\lfoot\ifnum\count0=\startpage\copyright\ \gtp\hfill\else
\agt, Volume \thevolumenumber\ (\thevolumeyear)\hfill\fi}
\def\@evenfoot{\@oddfoot}
\makeatother
\fi
\let\maketitlepage\makeagttitle
\let\makeshorttitle\maketitlepage
\let\maketitle\maketitlepage

\newwrite\gtoutfile
\long\gdef\makeheadfile{  
{\def\\{, }\def\s{ }
\immediate\openout\gtoutfile head.xxx
\immediate\write\gtoutfile{To: math@arxiv.org}
\immediate\write\gtoutfile{Subject: put OR rep NNNNN:ppppp}
\immediate\write\gtoutfile{--text follows this line--}
\immediate\write\gtoutfile{Proxy-for: \ifx\theasciiauthors\relax
\theauthors\else\theasciiauthors\fi\s<\ifx\theasciiemail\relax\theemail\else\theasciiemail\fi>}
\immediate\write\gtoutfile{\noexpand\\}
\immediate\write\gtoutfile{Authors: \ifx\theasciiauthors\relax
\theauthors\else\theasciiauthors\fi}
{\def\\{ }\immediate\write\gtoutfile{Title: \ifx\theasciititle\relax
\thetitle\else\theasciititle\fi}}
\immediate\write\gtoutfile{Subj-class: GT or SG, GR etc}
\immediate\write\gtoutfile{MSC-class: \theprimaryclass\ifx\thesecondaryclass\relax\else, \thesecondaryclass\fi}
\immediate\write\gtoutfile{Journal-ref: Algebr. Geom. Topol. \thevolumenumber\s
(\thevolumeyear) \startpage-\finishpage}
\immediate\write\gtoutfile{Comments: Published by Algebraic and
Geometric Topology at}
\immediate\write\gtoutfile{\s\s\s  http://www.maths.warwick.ac.uk/agt/AGTVol\thevolumenumber/agt-\thevolumenumber-\thepapernumber.abs.html}
\immediate\write\gtoutfile{\noexpand\\}
\immediate\write\gtoutfile{}
\ifx\theasciiabstract\relax
\immediate\write\gtoutfile{\theabstract}\else
\immediate\write\gtoutfile{\theasciiabstract}\fi
\immediate\write\gtoutfile{}
\immediate\write\gtoutfile{\noexpand\\}
\immediate\write\gtoutfile{}
\immediate\closeout\gtoutfile}}  

\def\maketitlepage{\makeagttitle\makeheadfile}
\let\makeshorttitle\maketitlepage
\let\maketitle\maketitlepage

\def\ifplaintex{\expandafter\ifx\csname documentclass\endcsname\relax}

\def\gtp{{\mathsurround=0pt\it $\cal G\mskip-2mu$eometry \&\ 
$\cal T\!\!$opology $\cal P\!$ublications}}  

\def\Addressesr{\bigskip
{\small \parskip 0pt \leftskip 0pt \rightskip 0pt plus 1fil \def\\{\par}
\sl\theaddress\par
\medskip
\rm Email:\stdspace\tt\theemail\hfill\rm Received:\qua\receiveddate \par}}

\def\recd{{\small Received:\qua\receiveddate\ifx\reviseddate\relax
\else\qquad Revised:\qua\reviseddate\fi\par}} 

\def\lognumber#1{\def\thelognumber{#1}}
\def\volumenumber#1{\def\thevolumenumber{#1}}
\def\volumeyear#1{\def\thevolumeyear{#1}}
\def\papernumber#1{\def\thepapernumber{#1}}
\def\pagenumbers#1#2{\def\startpage{#1}\def\finishpage{#2}}
\def\published#1{\def\publishdate{#1}}

\def\received#1{\def\receiveddate{#1}}

\def\accepted#1{\def\accepteddate{#1}}
\def\asciititle#1{\def\theasciititle{#1}}

\long\def\asciiabstract#1{\long\def\theasciiabstract{#1}}
\def\asciikeywords#1{\def\theasciikeywords{#1}}

\let\\\par\let\thelognumber\relax\let\thevolumenumber\relax
\let\thepapernumber\relax\let\thevolumeyear\relax\let\startpage\relax
\let\finishpage\relax\let\publishdate\relax\let\receiveddate\relax
\let\reviseddate\relax\let\accepteddate\relax\let\theasciititle\relax
\let\theasciiauthors\relax
\let\theasciiabstract\relax\let\theasciikeywords\relax

\let\theasciiemail\relax

\ifplaintex
\font\logobig=cmssbx10 scaled 3836
\font\logomed=cmssbx10 scaled 2557
\else
\font\logobig=cmssbx10 scaled 4200
\font\logomed=cmssbx10 scaled 2800
\fi

\long\def\makeagttitle{   
\count0=\startpage
\agt\hfill      
\hbox to 45truept{\vbox to 0pt{\vglue -13truept{\logomed A\kern -.37em{\logobig 
T}\kern -.38em G}\vss}\hss}
\break
{\small Volume \thevolumenumber\ (\thevolumeyear)
\startpage--\finishpage\nl
Published: \publishdate}

\vglue .25truein

{\parskip=0pt\leftskip 0pt plus
1fil\def\\{\par\smallskip}{\Large\bf\thetitle}\par\medskip} \vglue
0.05truein

{\parskip=0pt\leftskip 0pt plus 1fil\def\\{\par}{\sc\theauthors}
\par\medskip}%
 
\vglue 0.03truein 

{\small\leftskip 25truept\rightskip 25truept{\bf Abstract}\stdspace\theabstract

{\bf AMS Classification}\stdspace\theprimaryclass
\ifx\thesecondaryclass\relax\else; \thesecondaryclass\fi\par
{\bf Keywords}\stdspace \thekeywords\par}\vglue 7truept

}   

\ifplaintex
\font\phead=cmsl9 scaled 950
\font\pnum=cmbx10 scaled 913
\font\pfoot=cmsl9 scaled 950
\headline{\vbox to 0pt{\vskip -4.5mm\line{\small\phead\ifnum
\count0=\startpage ISSN 1472-2739 (on-line) 1472-2747 (printed)
\hfill {\pnum\folio}\else\ifodd\count0\def\\{ }%
\ifx\theshorttitle\relax\thetitle\else\theshorttitle\fi\hfill{\pnum\folio}
\else\def\\{ and }{\pnum\folio}\hfill\ifx\theshortauthors\relax\theauthors
\else\theshortauthors\fi\fi\fi}\vss}}
\footline{\vbox to 0pt{\vglue 0mm\line{\small\pfoot\ifnum\count0=\startpage
\copyright\ \gtp\hfill\else
\agt, Volume \thevolumenumber\ (\thevolumeyear)\hfill\fi}\vss}}
\else
\font=cmsl9 scaled 1050
\font\lnum=cmbx10 
\font=cmsl9 scaled 1050
\makeatletter
\def\@oddhead{{\small\lhead\ifnum\count0=\startpage ISSN 1472-2739 
(on-line) 1472-2747 (printed)\hfill {\lnum\number\count0}\else\ifodd\count0
\def\\{ }\ifx\theshorttitle\relax \thetitle \else\theshorttitle\fi\hfill
{\lnum\number\count0}\else\def\\{ and }{\lnum\number\count0}
\hfill\ifx\theshortauthors\relax 
\theauthors\else\theshortauthors\fi\fi\fi}}\def\@evenhead{\@oddhead}
\def\@oddfoot{\small\lfoot\ifnum\count0=\startpage\copyright\ \gtp\hfill\else
\agt, Volume \thevolumenumber\ (\thevolumeyear)\hfill\fi}
\def\@evenfoot{\@oddfoot}
\makeatother
\fi
\let\maketitlepage\makeagttitle
\let\makeshorttitle\maketitlepage
\let\maketitle\maketitlepage

\newwrite\gtoutfile
\long\gdef\makeheadfile{  
{\def\\{, }\def\s{ }
\immediate\openout\gtoutfile head.xxx
\immediate\write\gtoutfile{To: math@arxiv.org}
\immediate\write\gtoutfile{Subject: put OR rep NNNNN:ppppp}
\immediate\write\gtoutfile{--text follows this line--}
\immediate\write\gtoutfile{Proxy-for: \ifx\theasciiauthors\relax
\theauthors\else\theasciiauthors\fi\s<\ifx\theasciiemail\relax\theemail\else\theasciiemail\fi>}
\immediate\write\gtoutfile{\noexpand\\}
\immediate\write\gtoutfile{Authors: \ifx\theasciiauthors\relax
\theauthors\else\theasciiauthors\fi}
{\def\\{ }\immediate\write\gtoutfile{Title: \ifx\theasciititle\relax
\thetitle\else\theasciititle\fi}}
\immediate\write\gtoutfile{Subj-class: GT or SG, GR etc}
\immediate\write\gtoutfile{MSC-class: \theprimaryclass\ifx\thesecondaryclass\relax\else, \thesecondaryclass\fi}
\immediate\write\gtoutfile{Journal-ref: Algebr. Geom. Topol. \thevolumenumber\s
(\thevolumeyear) \startpage-\finishpage}
\immediate\write\gtoutfile{Comments: Published by Algebraic and
Geometric Topology at}
\immediate\write\gtoutfile{\s\s\s  http://www.maths.warwick.ac.uk/agt/AGTVol\thevolumenumber/agt-\thevolumenumber-\thepapernumber.abs.html}
\immediate\write\gtoutfile{\noexpand\\}
\immediate\write\gtoutfile{}
\ifx\theasciiabstract\relax
\immediate\write\gtoutfile{\theabstract}\else
\immediate\write\gtoutfile{\theasciiabstract}\fi
\immediate\write\gtoutfile{}
\immediate\write\gtoutfile{\noexpand\\}
\immediate\write\gtoutfile{}
\immediate\closeout\gtoutfile}}  

\def\maketitlepage{\makeagttitle\makeheadfile}
\let\makeshorttitle\maketitlepage
\let\maketitle\maketitlepage